\documentclass[11pt,twoside]{amsart}
\usepackage{amsfonts,amssymb} 
\renewcommand{\baselinestretch}{1.15}
\def\botrel#1#2{\mathrel{\mathop{#1}\limits_{#2}}}

\def\sqr#1#2{{\vcenter{\hrule height.#2pt \hbox{\vrule width.#2pt
height#1pt \kern#1pt \vrule width.#2pt} \hrule height.#2pt}}}
\def\square{\mathchoice\sqr64\sqr64\sqr{4}3\sqr{3}3}
\def\qed{\hfill$\square$}
\def\demo{\noindent{\bf Proof: }}

\newtheorem{Proposition}{\bf Proposition}[section]
\newtheorem{Theorem}[Proposition]{\bf Theorem}
\newtheorem{Lemma}[Proposition]{\bf Lemma}
\newtheorem{Corollary}[Proposition]{\bf Corollary}

\newtheorem{Remark}[Proposition]{\bf Remark}
\newtheorem{Example}[Proposition]{\bf Example}

\newtheorem*{Theorem/id=rn}{\bf Theorem \ref{id=rn}}
\newtheorem*{Theorem/modules}{\bf Theorem \ref{modules}}
\newtheorem*{Theorem/ideals}{\bf Theorem \ref{ideals}}
\newtheorem*{Theorem/ideals2}{\bf Theorem \ref{ideals2}}

\newcommand{\id}{\mbox{$\mathrm{id}$}}
\newcommand{\rt}{\mbox{$\mathrm{rt}$}}
\newcommand{\rn}{\mbox{$\mathrm{rn}$}}
\newcommand{\ud}{\mbox{$\mathrm{d}$}}
\newcommand{\reg}{\mbox{$\mathrm{reg}$}}
\newcommand{\rees}[1]{\mbox{$\mathcal{R}(#1)$}}
\newcommand{\reesq}[2]{\mbox{$\mathcal{R}_{#1}(#2)$}}
\newcommand{\reesm}[2]{\mbox{$\mathcal{R}(#1;#2)$}}
\newcommand{\reesqm}[3]{\mbox{$\mathcal{R}_{#1}(#2;#3)$}}
\newcommand{\agr}[1]{\mbox{$\mathcal{G}(#1)$}}
\newcommand{\agrm}[2]{\mbox{$\mathcal{G}(#1;#2)$}}
\newcommand{\fcone}[1]{\mbox{$\mathcal{F}_{\mathfrak{m}}(#1)$}}
\newcommand{\fconem}[2]{\mbox{$\mathcal{F}_{\mathfrak{m}}(#1;#2)$}}
\newcommand{\fa}{\mbox{$\mathfrak{a}$}}
\newcommand{\fm}{\mbox{$\mathfrak{m}$}}
\newcommand{\fn}{\mbox{$\mathfrak{n}$}}
\newcommand{\fp}{\mbox{$\mathfrak{p}$}}

\newcommand{\bq}{\mbox{$\mathbb{Q}$}}

\begin{document}

\title{Integral degree of a ring and reduction numbers}

\author{{\sc Jos\'e M. Giral and Francesc Planas-Vilanova}}

\date{\today}

\begin{abstract}
The supremum of reduction numbers of ideals having principal
reductions is expressed in terms of the integral degree, a new
invariant of the ring, which is finite provided the ring has finite
integral closure. As a consequence, one obtains bounds for the
Castelnuovo-Mumford regularity of the Rees algebra and for the
Artin-Rees numbers.
\end{abstract}

\maketitle

\section{Introduction}

Let $A$ be a commutative noetherian ring with identity, let $I$ be an
ideal of $A$ and let $M$ be a finitely generated $A$-module. An ideal
$J\subset I$ is said to be a reduction of $I$ with respect to $M$ if
$I^{n+1}M=JI^{n}M$ for some integer $n\geq 0$. The least such integer
$n\geq 0$ is called the $J$-{\em reduction number of} $I$ {\em with
respect to} $M$ and is denoted by ${\rm rn}_{J}(I;M)$. If $M=A$, the
phrase ``with respect to $M$'' is omitted and one writes ${\rm
rn}_{J}(I)$. Clearly, if $J$ is a reduction of $I$, then $J$ is a
reduction of $I$ with respect to $M$ and ${\rm rn}_{J}(I;M)\leq {\rm
rn}_{J}(I)$. If $I$ is regular (i.e. $I$ contains a non zero divisor)
and $J$ is a principal reduction, then $\rn_{J}(I)$ is independent of
the given principal reduction and is denoted by $\rn(I)$.

In the last decade there has been a great deal of attention to finding
bounds on the reduction number (see e.g. \cite{aht}, \cite{cpv},
\cite{dgh2}, \cite{dkv}, \cite{hoa} \cite{rossi}, \cite{vasconcelos1},
\cite{vasconcelos2}, \cite{vasconcelos3}, by no means a complete list
of references). If the ideals have principal reductions and the
integral closure is finite, d'Anna, Guerrieri and Heinzer gave an
absolute bound for the reduction number in terms of the minimal number
of generators of the integral closure (\cite{dgh2}, Corollary 5.2). In
this paper we express the supremum of reduction numbers of ideals
having principal reductions in terms of the following new invariant
associated to the ring $A$, provided $A$ contains the field of
rational numbers $\bq$. If $A\subset B$ is a ring extension and $b\in
B$ is integral over $A$, let the {\em integral degree} of $b$ over $A$
be
\begin{eqnarray*}
{\rm id}_{A}(b)={\rm min}\{ n\geq 1\mid b\mbox{ satisfies an integral
equation of degree }n\}.
\end{eqnarray*}
If $A\subset B$ is an integral extension, the {\em integral degree} of
$B$ over $A$ is defined as
\begin{eqnarray*}
\ud _{A}(B)={\rm sup}\{{\rm id}_{A}(b)\mid b\in B\},
\end{eqnarray*}
When $B$ is taken to be $\overline{A}$, the integral closure of $A$ in
its total quotient ring, $\ud_{A}(\overline{A})$ is just called the
{\em integral degree} of $A$. We will prove that if $A$ has finite
integral closure then it has also finite integral degree and that the
converse is not true in general. Our main result is:

\begin{Theorem/id=rn}
Let $A$ be a noetherian ring, $A\supset\bq$. Then
\begin{eqnarray*}
\ud_{A}(\overline{A})={\rm sup}\, \{ {\rm rn}(I)\mid I\mbox{ regular
ideal of }A\mbox{ having a principal reduction}\} +1.
\end{eqnarray*}
\end{Theorem/id=rn}

In Theorem~\ref{id=rn} it is possible to replace $\rn(I)+1$ by either
$\reg(\rees{I})+1$ or $\rt(I)$, $\reg(\rees{I})$ and $\rt(I)$ being
the Castel\-nuovo-Mumford regularity of the Rees algebra of $I$ and
the relation type of $I$, respectively.

It is known that Artin-Rees numbers are bounded by the relation type
and that, in some particular cases, the relation type can be bounded
by reduction numbers. Having in mind this idea and as a consequence of
Theorem~\ref{id=rn}, we get the following results in the context of
uniform Artin-Rees properties.

\begin{Theorem/modules}
Let $A$ be a noetherian ring with finite integral degree $\ud
_{A}(\overline{A})=d$. Suppose that $A\supset\bq$. Let $N\subset M$ be
two finitely generated $A$-modules. Let $I$ be a regular ideal of $A$
having a principal reduction generated by a $d$-sequence with respect
to $M/N$. Then, for every integer $n\geq d$,
\begin{eqnarray*}
I^{n}M\cap N=I^{n-d}(I^{d}M\cap N).
\end{eqnarray*}
\end{Theorem/modules}
  
In other words, $\ud_{A}(\overline{A})$ is a uniform Artin-Rees number
for the pair $N\subset M$ and the whole set of regular ideals having
principal reductions generated by a $d$-sequence with respect to
$M/N$.  Our ideal-theoretic version is the following.

\begin{Theorem/ideals}
Let $A$ be a noetherian ring, $A\supset\bq$. Let $\fa$ be an ideal of
$A$ such that $A/\fa$ has finite integral degree
$\ud_{A/\fa}(\overline{A/\fa})=d$. Let $I$ be an ideal of $A$ such
that $IA/\fa$ has an $A/\fa$-regular principal reduction. Then, for
every integer $n\geq d$,
\begin{eqnarray*}
I^{n}\cap \fa=I^{n-d}(I^{d}\cap \fa).
\end{eqnarray*}
\end{Theorem/ideals}

If the integral degree of $\overline{A/\fa}$ is not finite or if
$IA/\fa$ has no principal reduction, then there may not exist such a
uniform Artin-Rees number (Example~\ref{ehe} and
Example~\ref{wange}). This will be seen by using an example of
Eisenbud and Hochster in \cite{eh}, the work where they raised the
uniform Artin-Rees conjecture, and an example of Wang in \cite{wang}
(see also \cite{ocarroll1}, \cite{do}, \cite{ocarroll2},
\cite{huneke}, \cite{planas} for more information). On the other hand,
it is well known that there exists a uniform Artin-Rees number for the
set of principal ideals of a noetherian ring and that, in general,
there does not exist a uniform Artin-Rees number for the set of three
generated ideals (see the work of O'Carroll in \cite{ocarroll2} and
the just mentioned example of Wang in \cite{wang}). Therefore, it
remained to study if there exists a uniform Artin-Rees number for the
whole set of two-generated ideals (without any other assumption on the
ideals). We obtain a slightly weaker uniform Artin-Rees property for
the set of two-generated regular ideals (not true anymore for the set
of three-generated ideals, see Example \ref{wange}). Concretely,

\begin{Theorem/ideals2}
Let $(A,\fm)$ be a noetherian local ring with infinite residue
field. Let $\fa$ be an ideal of $A$ such that $A/\fa$ has finite
integral degree $\ud_{A/\fa}(\overline{A/\fa})=d$. Let $I$ be a
two-generated ideal of $A$ such that $IA/\fa$ is
$A/\fa$-regular. Then, for every $n\geq d$,
\begin{eqnarray*}
I^{n}\cap \fa =I^{n-d}(I^{d}\cap \fa) +\fm I^{n}\cap \fa .
\end{eqnarray*}
\end{Theorem/ideals2}

The paper is organized as follows. Sections~\ref{armj}, \ref{rtmj},
\ref{cmr} and \ref{rtrn} are devoted to the following four invariants
and the relationship among them: Artin-Rees numbers modulo an ideal,
relation type of a standard module, Castelnuovo-Mumford regularity and
reduction number with respect to a module. Concretely, in Section
\ref{armj} we introduce the Artin-Rees number $s_{J}(N,M;I)$ of an
ideal $I$, two finitely generated $A$-modules $N\subset M$ and modulo
another ideal $J$. This number is the minimum integer $s\geq 0$ such
that $I^{n}M\cap N=I^{n-s}(I^{s}M\cap N)+JI^{n}M\cap N$ for all $n\geq
s+1$, and thus it controls the weaker Artin-Rees property of
Theorem~\ref{ideals2}. Following the ideas in \cite{planas}, in
Section~\ref{rtmj} we bound above the Artin-Rees number $s_{J}(N,M;I)$
by the relation type of the Rees module $\reesqm{J}{I}{M/N}=(\oplus
_{n\geq 0}I^{n}M/N)\otimes A/J$. Section~\ref{cmr} is dedicated to
recalling some definitions around Castelnuovo-Mumford regularity and
formulating an extension to modules of some results of Trung in
\cite{trung1} and \cite{trung2}. In Section~\ref{rtrn}, we prove that
the relation type of an ideal $I$ with respect to a module $M$, ${\rm
rt}(I;M)$, is bounded above by $\rn_{J}(I;M)+\rt(I;J^{r}M)$, where $J$
is a reduction of $I$ with respect to $M$ and $r:=\rn_{J}(I;M)$ is the
$J$-reduction number of $I$ with respect to $M$. If $J$ is generated
by a complete $d$-sequence with respect to $I$ and $M$ (the
terminology is explained in Section~\ref{rtrn}), then the relation
type of $J$ with respect to $I^{r}M$ satisfies ${\rm rt}(J;I^{r}M)=1$
and $\reg(\reesm{I}{M})=\rn_{J}(I;M)$. Thus one has the inequality
${\rm rt}(I;M)\leq {\rm rn}_{J}(I;M)+1$, which is well-known for the
case $M=A$ and $J$ a principal reduction of a regular ideal $I$ (see
the work of d'Anna, Guerrieri and Heinzer~\cite{dgh1}, Huckaba
\cite{huckaba2}, \cite{huckaba3}, Schenzel \cite{schenzel} and Trung
\cite{trung1},\cite{trung2}). In Section \ref{id}, we introduce and
study $\ud_{A}(\overline{A})$, the integral degree of $A$, a new
invariant associated to the ring $A$. We prove that if $A$ has finite
integral closure then it also has finite integral degree. The
ingenious example of Akizuki (\cite{akizuki}, see also \cite{reid})
provides us with an example of a one-dimensional noetherian local
domain $A$ with finite integral degree but infinite integral
closure. In Section~\ref{idrn} we prove the main result of the paper,
namely, $\ud_{A}(\overline{A})$ is equal to the supremum of the
reduction numbers plus one (or else the Castelnuovo-Mumford regularity
of the Rees algebra plus one or the relation type) of regular ideals
having principal reductions.  Finally, in Section \ref{uarp} we prove
all the results concerning Artin-Rees numbers.

All rings will be commutative and with identity. As usual, $M$-regular
will mean not contained in the set of zero divisors of $M$ and $\mu$
will stand for the minimal number of generators.

\section{Artin-Rees modulo an ideal}\label{armj}

Let us introduce a slight variant of the Artin-Rees Lemma which will
be very useful. Let $A$ be a noetherian ring, $I$ an ideal of $A$ and
$N\subseteq M$ two finitely generated $A$-modules. The Artin-Rees
Lemma assures that there exists an integer $s\geq 0$, depending on
$N$, $M$ and $I$, such that for all $n\geq s$,
\begin{eqnarray*}
I^{n}M\cap N=I^{n-s}(I^{s}M\cap N).
\end{eqnarray*}
In particular, for any ideal $J$ of $A$, one obtains what we will call
{\em Artin-Rees modulo} $J$:
\begin{eqnarray*}
I^{n}M\cap N=I^{n-s}(I^{s}M\cap N)+JI^{n}M\cap N .
\end{eqnarray*}
For every integer $n\geq 1$, let 
\begin{eqnarray*}
E_{J}(N,M;I)_{n}=\frac{I^{n}M\cap N}{I(I^{n-1}M\cap
N)+JI^{n}M\cap N}\, .
\end{eqnarray*}
For easy reference, and without proof, we state the Artin-Rees lemma
modulo $J$.

\begin{Lemma}\label{sj} 
Let $A$ be a ring, $I,J$ ideals of $A$ and $N\subseteq M$ two
$A$-modules.  Set 
\begin{eqnarray*}
s_{J}(N,M;I)={\rm min}\{ s\geq 0\mid
E_{J}(N,M;I)_{n}=0 \mbox{ for all } n\geq s+1\}.
\end{eqnarray*}
Then, the following conditions are equivalent:
\begin{itemize}
\item[$(i)$] $I^{n}M\cap N=I^{n-s}(I^{s}M\cap N)+JI^{n}M\cap N$ for
all $n\geq s+1$.
\item[$(ii)$] $s_{J}(N,M;I)\leq s$.
\end{itemize}
If $A$ is noetherian and $N\subseteq M$ are finitely-generated
$A$-modules, then $s_{J}(N,M;I)$ is finite.
\end{Lemma}

If $J=0$, we recover the standard notion of Artin-Rees and simply
write $s(N,M;I)$. Remark that if $J_{1}\subset J_{2}$ are two ideals,
and for $n\geq 1$, there is a natural epimorphism
$E_{J_{1}}(N,M;I)_{n}\to E_{J_{2}}(N,M;I)_{n}$ and thus
$s_{J_{2}}(N,M;I)\leq s_{J_{1}}(N,M;I)\leq s(N,M;I)$.

\begin{Remark}\label{sjs0}{\rm 
If $A$ is noetherian, $J\subset I$ are two ideals of $A$ contained in
the Jacobson radical of $A$, and $N\subset M$ are two finitely
generated $A$-modules, then $s_{J}(N,M;I)=s(N,M;I)$.  }
\end{Remark}

\demo Since $0\subset J\subset I$, then $s_{I}(N,M;I)\leq
s_{J}(N,M;I)\leq s(N,M;I)$. It is enough to see $s(N,M;I)\leq
s_{I}(N,M;I)$. Set $s=s_{I}(N,M;I)$, so $I^{n}M\cap
N=I^{n-s}(I^{s}M\cap N)+I^{n+1}M\cap N$ for all $n\geq s+1$. Then
$I^{n+1}M\cap N=I^{n+1-s}(I^{s}M\cap N)+I^{n+2}M\cap N$ and
substituting the second equality in the first, $I^{n}M\cap
N=I^{n-s}(I^{s}M\cap N)+I^{n+1-s}(I^{s}M\cap N)+I^{n+2}M\cap
N=I^{n-s}(I^{s}M\cap N)+I^{n+2}M\cap N$. Inductively, $I^{n}M\cap
N=\cap _{k\geq 1}(I^{n-s}(I^{s}M\cap N)+I^{n+k}M\cap N)\subset \cap
_{k\geq 1}(P+I^{n+k}M)$, where $P=I^{n-s}(I^{s}M\cap N)\subset
I^{n}M\cap N\subset M$. But,
\begin{eqnarray*}
\frac{\bigcap _{k\geq 1}(P+I^{n+k}M)}{P}= \bigcap _{k\geq
  1}\left(\frac{P+I^{n+k}M}{P}\right)=\bigcap _{k\geq 1}I^{n+k}(M/P),
\end{eqnarray*}
which is zero by Krull's intersection theorem. Therefore, $\cap
_{k\geq 1}(P+I^{n+k}M)=P$ and $I^{n}M\cap N=I^{n-s}(I^{s}M\cap N)$. \qed

\section{Relation type modulo an ideal}\label{rtmj}

A standard $A$-algebra is a commutative graded algebra $U=\oplus
_{n\geq 0}U_{n}$, with $U_{0}=A$ and $U$ generated by the elements of
degree 1. The Rees algebra of $I$ is the standard $A$-algebra
$\rees{I}=\oplus _{n\geq 0}I^{n}$. For any ideal $J$ of $A$, the {\em
Rees algebra of} $I$ {\em modulo} $J$ will be the standard
$A/J$-algebra $\reesq{J}{I}=\rees{I}\otimes A/J=\oplus _{n\geq 0}
I^{n}/JI^{n}$. Taking $J=I$, we recover the associated graded ring of
$I$, $\reesq{I}{I}=\agr{I}=\oplus _{n\geq 0}I^{n}/I^{n+1}$, and taking
$J=\fm$ a maximal ideal of $A$, we recover the fiber cone of $I$,
$\reesq{\fm}{I}=\fcone{I}=\oplus _{n\geq 0}I^{n}/\fm I^{n}$.

A standard $U$-module will be a graded $U$-module $F=\oplus _{n\geq
0}F_{n}$ such that $F_{n}=U_{n}F_{0}$ for all $n\geq 0$. The Rees
module of $I$ with respect to $M$ is the standard $\rees{I}$-module
$\reesm{I}{M}=\oplus _{n\geq 0}I^{n}M$. For any ideal $J$ of $A$, the
{\em Rees module of} $I$ {\em with respect to} $M$ {\em and modulo}
$J$ will be the standard $\reesq{J}{I}$-module
$\reesqm{J}{I}{M}=\reesm{I}{M}\otimes A/J=\oplus _{n\geq 0}
I^{n}M/JI^{n}M$. Taking $J=I$, we recover the associated graded module
of $I$ with respect to $M$, $\reesqm{I}{I}{M}=\agrm{I}{M}=\oplus
_{n\geq 0}I^{n}M/I^{n+1}M$ and taking $J=\fm$ a maximal ideal of $A$,
we recover the fiber cone of $I$ with respect to $M$,
$\reesqm{\fm}{I}{M}=\fconem{I}{M}=\oplus _{n\geq 0}I^{n}M/\fm I^{n}M$.

Given two standard $U$-modules $F$, $G$ and $\varphi :G\rightarrow F$,
a surjective graded morphism of $U$-modules, put $E(\varphi )_{n}=
{\rm ker}\varphi _{n}/U_{1}{\rm ker}\varphi _{n-1}$ for $n\geq
2$. Consider $\gamma :{\bf S}(U_{1})\otimes F_{0}\buildrel \alpha
\otimes 1\over \rightarrow U\otimes F_{0}\rightarrow F$, where $\alpha
:{\bf S}(U_{1})\rightarrow U$ is the canonical symmetric presentation
of $U$ and $U\otimes F_{0}\rightarrow F$ is the structural
morphism. For $n\geq 2$, the {\em module of effective $n$-relations
of} $F$ is $E(F)_{n}=E(\gamma )_{n}= {\rm ker}\gamma _{n}/U_{1}{\rm
ker}\gamma _{n-1}$. The {\em relation type of} $F$ is ${\rm
rt}(F)={\rm min}\{ r\geq 1\mid E(F)_{n}=0\mbox{ for all }n\geq r+1\}$,
which is finite if $A$ is noetherian, $U$ is a finitely generated
algebra and $F$ is a finitely generated $U$-module. It can be shown
that the module of effective $n$-relations, $n\geq 2$, and the
relation type do not depend on the chosen symmetric presentation
(\cite{planas}, Definition 2.4, see also \cite{dgh1},
\cite{vasconcelos3}, \cite{wang}). In particular, in order to find the
effective relations of $F$ and its relation type, one can always take
a presentation of $F$ as a quotient of a polynomial module with
coeficients in $F_{0}$.

The {\em module of effective} $n$-{\em relations of} $I$ {\em with
respect to} $M$ is $E(I;M)_{n}=E(\reesm{I}{M})_{n}$ and the relation
type of $I$ with respect to $M$ is ${\rm rt}(I;M)={\rm
rt}(\reesm{I}{M})$. For any ideal $J$ of $A$, the {\em module of
effective $n$-relations of} $I$ {\em with respect to} $M$ {\em and
modulo} $J$ will be $E_{J}(I;M)_{n}=E(\reesqm{J}{I}{M})_{n}$ and the
relation type of $I$ with respect to $M$ and modulo $J$ will be ${\rm
rt}_{J}(I;M)={\rm rt}(\reesqm{J}{I}{M})$.  If $M=A$, then we omit the
phrase ``with respect to $M$'' and simply write $E(I)_{n}$, ${\rm
rt}(I)$, $E_{J}(I)_{n}$ and ${\rm rt}_{J}(I)$.

\begin{Remark}\label{coef-mod}{\rm 
Let $A$ be a ring, $J,I,\fa$, ideals of $A$ and $M$ an $A$-module.
\begin{itemize}
\item[(1)] Then ${\rm rt}(I;A/\fa)={\rm rt}(IA/\fa)={\rm
rt}(IA/\fa;A/\fa)$.
\item[(2)] ${\rm rt}_{J}(I;M)\leq {\rm rt}(I;M)$.
\item[(3)] If $A$ is noetherian, $J\subset I$ and $M$ is finitely
generated, then ${\rm rt}_{J}(I;M)={\rm rt}(I;M)$.
\end{itemize}
}\end{Remark}

\demo $\reesm{I}{A/\fa}=\rees{IA/\fa}=\reesm{IA/\fa}{A/\fa}$.
Moreover, the relation type of the standard $\rees{I}$-module
$\reesm{I}{A/\fa}$, the relation type of the standard $A/\fa$-algebra
$\rees{IA/\fa}$ and the relation type of the standard $A/\fa$-module
$\reesm{IA/\fa}{A/\fa}$ all coincide (see \cite{planas}, Remark
2.5). This proves $(1)$. The proof of $(2)$ and $(3)$ follow from
\cite{planas}, Remark 2.7 (and in contrast to Remark \ref{sjs0}, here
we do not need $I$ to be included in the Jacobson radical). \qed

\medskip

Next we show the relation between $E(I;M)_{n}$ and $E_{J}(I;M)_{n}$
and describe these modules for the two-generated regular case.

\begin{Proposition}\label{E(x,y)}
Let $A$ be a ring, $I$ and $J$ ideals of $A$ and $M$ an
$A$-module. Then, for every integer $n\geq 2$, there exists an exact
sequence of $A$-modules:
\begin{eqnarray*}
E(I;JM)_{n}\longrightarrow E(I;M)_{n}\longrightarrow
E_{J}(I;M)_{n}\rightarrow 0.
\end{eqnarray*}
In particular, if $I=(x,y)$ is two-generated and $x$ is $M$-regular,
then, for every $n\geq 2$,
\begin{eqnarray*}
E_{J}(I;M)_{n}=\frac{(xI^{n-1}M:y^{n})}{ (xI^{n-1}JM:y^{n})\cap
JM+(xI^{n-2}M:y^{n-1})}.
\end{eqnarray*}
\end{Proposition}

\demo Let $f:P\rightarrow I$ be a presentation of $I$, with $P$ a free
$A$-module, and, for every $n\geq 2$, consider the following
commutative diagram:

\begin{picture}(330,125)(0,0)

\put(75,100){\makebox(0,0){\mbox{\footnotesize
$\Lambda_{2}(P)\otimes I^{n-2}JM$}}}
\put(210,100){\makebox(0,0){\mbox{\footnotesize
$P\otimes I^{n-1}JM$}}}
\put(310,100){\makebox(0,0){\mbox{\footnotesize
$I^{n}JM$}}}
\put(370,100){\makebox(0,0){$0$}}

\put(130,100){\vector(1,0){40}}
\put(250,100){\vector(1,0){40}}
\put(330,100){\vector(1,0){30}}

\put(75,60){\makebox(0,0){\mbox{\footnotesize
$\Lambda_{2}(P)\otimes I^{n-2}M$}}}
\put(210,60){\makebox(0,0){\mbox{\footnotesize
$P\otimes I^{n-1}M$}}}
\put(310,60){\makebox(0,0){\mbox{\footnotesize
$I^{n}M$}}}
\put(370,60){\makebox(0,0){$0$}}

\put(130,60){\vector(1,0){40}}
\put(250,60){\vector(1,0){40}}
\put(330,60){\vector(1,0){30}}

\put(75,20){\makebox(0,0){\mbox{\footnotesize 
$\Lambda_{2}(P)\otimes I^{n-2}M/I^{n-2}JM$}}}
\put(210,20){\makebox(0,0){{\footnotesize 
$P\otimes I^{n-1}M/I^{n-1}JM$}}}
\put(310,20){\makebox(0,0){\mbox{\footnotesize 
$I^{n}M/I^{n}JM$}}}
\put(370,20){\makebox(0,0){$0$}}

\put(140,20){\vector(1,0){20}}
\put(260,20){\vector(1,0){20}}
\put(340,20){\vector(1,0){20}}

\put(75,50){\vector(0,-1){20}}
\put(75,50){\vector(0,-1){16}}
\put(210,50){\vector(0,-1){20}}
\put(210,50){\vector(0,-1){16}}
\put(310,50){\vector(0,-1){20}}
\put(310,50){\vector(0,-1){16}}
\put(75,90){\vector(0,-1){20}}
\put(210,90){\vector(0,-1){20}}
\put(310,88){\makebox(0,0){{\footnotesize $\lor$}}}
\put(310,86){\vector(0,-1){16}}

\put(150,110){\makebox(0,0){\mbox{\footnotesize
$\partial^{\prime}_{2,n}$}}}
\put(150,70){\makebox(0,0){\mbox{\footnotesize $\partial _{2,n}$}}}
\put(150,30){\makebox(0,0){\mbox{\footnotesize
$\overline{\partial}_{2,n}$}}}
\put(270,110){\makebox(0,0){\mbox{\footnotesize
$\partial^{\prime}_{1,n}$}}}
\put(270,70){\makebox(0,0){\mbox{\footnotesize $\partial _{1,n}$}}}
\put(270,30){\makebox(0,0){\mbox{\footnotesize
$\overline{\partial}_{1,n}$}}}

\end{picture}

\noindent The top, middle and bottom rows of these diagrams represent
the last three nonzero terms of the $n$-th homogeneous part of the
Koszul complexes induced by the ${\bf S}(P)$-linear forms $P\otimes
\reesm{I}{JM}\rightarrow \reesm{I}{JM}$, $P\otimes
\reesm{I}{M}\rightarrow \reesm{I}{M}$ and $P\otimes
\reesqm{J}{I}{M}\rightarrow \reesqm{J}{I}{M}$. The differentials are
defined as usual: $\partial_{2,n}((x\wedge y)\otimes z)=y\otimes
xz-x\otimes yz$ and $\partial_{1,n}(x\otimes t)=xt$, $x,y\in P$, $z\in
I^{n-2}M$, $t\in I^{n-1}M$; $\partial^{\prime}_{i,n}$ and
$\overline{\partial}_{i,n}$ are defined analogously (see
e.g. \cite{bh}, Definition~1.6.1). The vertical morphisms are induced
by the obvious inclusions and quotients and define morphisms of
complexes. By a similar reasoning to that in \cite{planas},
Proposition 2.6, the first homology groups of these complexes are,
respectively, ${\rm ker}\partial^{\prime}_{1,n}/{\rm
im}\partial^{\prime}_{2,n}=E(I;JM)_{n}$, ${\rm ker}\partial_{1,n}/{\rm
im}\partial_{2,n}=E(I;M)_{n}$ and ${\rm
ker}\overline{\partial}_{1,n}/{\rm
im}\overline{\partial}_{2,n}=E_{J}(I;M)_{n}$. The exact sequence we
seek is nothing else but the short exact sequence induced in homology.

If $I=(x,y)$ with $x$ an $M$-regular element, take $P=A^{2}$ and
$f:P\rightarrow I$ with $f(1,0)=x$ and $f(0,1)=y$. Then, the middle
row becomes isomorphic to the complex:
\begin{eqnarray*}
I^{n-2}M\buildrel \partial_{2,n}\over\longrightarrow I^{n-1}M\oplus
I^{n-1}M\buildrel\partial_{1,n}\over\longrightarrow I^{n}M\rightarrow 0,
\end{eqnarray*}
with differentials $\partial_{2,n}(u)=(-yu,xu)$ and
$\partial_{1,n}(z,t)=xz+yt$. Take $(z,t)=(\sum a_{i}u_{i},\sum
b_{i}v_{i})$, $(z,t)\in {\rm ker}\partial_{1,n}$, with $a_{i},b_{i}\in
I^{n-1}$, $u_{i},v_{i}\in M$ and $xz+yt=0$. Write
$b_{i}=c_{i}y^{n-1}+d_{i}x$, $c_{i}\in A$ and $d_{i}\in I^{n-2}$. Then
\begin{eqnarray*}
y^{n}\sum c_{i}v_{i}=y\sum b_{i}v_{i}-y\sum d_{i}xv_{i}=yt-x\sum
d_{i}yv_{i}=-x(z-\sum d_{i}yv_{i}).
\end{eqnarray*}
Thus $\sum c_{i}v_{i}\in (xI^{n-1}M:y^{n})$.
Consider 
\begin{eqnarray*}
\varphi :{\rm ker}\partial_{1,n}\longrightarrow
\frac{(xI^{n-1}M:y^{n})}{(xI^{n-2}M:y^{n-1})} \; ,
\end{eqnarray*}
defined by $\varphi (z,t)=\overline{\sum c_{i}v_{i}}$. It is not
difficult to see that $\varphi$ is well-defined, surjective and ${\rm
im}\partial_{2,n}\subset {\rm ker}\varphi$. Moreover, if $x$ is
$M$-regular, then ${\rm ker}\varphi\subset {\rm
im}\partial_{2,n}$. Thus
\begin{eqnarray*}
E(I;M)_{n}=\frac{(xI^{n-1}M:y^{n})}{(xI^{n-2}M:y^{n-1})}.
\end{eqnarray*}
Using the former exact sequence of modules of effective relations, one
deduces the expression of $E_{J}(I;M)_{n}$. \qed

\medskip

Next we compare the Artin-Rees number modulo $J$ with the relation
type modulo $J$.

\begin{Proposition}\label{slrt}
Let $A$ be a ring, $I$ and $J$ two ideals of $A$ and $N\subset M$ two
$A$-modules. Then
\begin{eqnarray*}
s_{J}(N,M;I)\leq {\rm rt}_{J}(I;M/N)\leq {\rm max}({\rm
rt}_{J}(I;M),s_{J}(N,M;I)).
\end{eqnarray*}
\end{Proposition}

\demo Take $F=\reesqm{J}{I}{M/N}$, $G=\reesqm{J}{I}{M}$ and $H={\bf
S}(I/JI)\otimes M$ and $\varphi :G\rightarrow F$ and $\gamma
:H\rightarrow G$ induced by the natural surjective graded morphisms
$\reesm{I}{M}\rightarrow \reesm{I}{M/N}$ and ${\bf S}(I)\otimes
M\rightarrow \reesm{I}{M}$. By \cite{planas}, Lemma 2.3, for every
integer $n\geq 2$, one has the short exact sequence of $A$-modules:
\begin{eqnarray*}
E(\gamma )_{n}\rightarrow E(\varphi\circ\gamma )_{n}\rightarrow
E(\varphi )_{n}\rightarrow 0.
\end{eqnarray*}
But $E(\gamma )_{n}=E_{J}(I;M)_{n}$ and $E(\varphi\circ\gamma
)_{n}=E_{J}(I;M/N)_{n}$ and a short computation shows that $E(\varphi
)_{n}=E_{J}(N,M;I)_{n}$. From the exact sequence we obtain the desired
inequalities. \qed

\section{Castelnuovo-Mumford regularity}\label{cmr}

The purpose of this section is to recall some definitions and
formulate, in order to use them subsequently, a generalization to
modules of some results of Trung in \cite{trung1} and
\cite{trung2}. Being natural extensions of his results, we omit or
just sketch the proofs. Let $A$ be a noetherian ring and $U=\oplus
_{n\geq 0}U_{n}$ a finitely generated standard $A$-algebra. Let
$F=\oplus _{n\geq 0}F_{n}$ be a standard $U$-module. Define
\begin{eqnarray*}
a(F)=\left\{ \begin{array}{ll} \mbox{\rm max}\{ n\geq 0\mid F_{n}\neq
0\}&\mbox{ if }F\neq 0.\\ -\infty &\mbox{ if }F=0.\end{array}\right.
\end{eqnarray*}
Let $U_{+}=\oplus_{n>0}U_{n}$ be the irrelevant ideal of $U$. If
$i\geq 0$, denote
by
\begin{eqnarray*}
a_{i}(F)=a(H^{i}_{U_{+}}(F)),
\end{eqnarray*}
where $H^{i}_{U_{+}}(\cdot )$ denotes the $i$-th local cohomology
functor with respect to the ideal $U_{+}$. The {\em
Castelnuovo-Mumford regularity of $F$} is defined to be
\begin{eqnarray*}
\reg{(F)}=\mbox{\rm max}\{ a_{i}(F)+i\mid i\geq 0\}
\end{eqnarray*}
(see e.g. \cite{bs}, 15.2.9, \cite{schenzel}, \cite{trung2}). We shall
mainly be concerned with the case $U=\rees{I}$, the Rees algebra of an
ideal $I$ of $A$, and $F=\reesm{I}{M}$, the Rees module of $I$ with
respect to a finitely generated $A$-module $M$. In particular, if
$M\neq 0$, then $\reg{(F)}\neq -\infty$ (see e.g. \cite{bs}, 15.2.13).

A sequence ${\bf z}=z_{1},\ldots ,z_{s}$ of homogeneous elements
of $U$ is called {\em $n$-regular with respect to $F$} if, for
all $i=1,\ldots, s$,
\begin{eqnarray*}
((z_{1},\ldots ,z_{i-1})F:z_{i})_{n}=((z_{1},\ldots ,z_{i-1})F)_{n}.
\end{eqnarray*}
The least integer $m\geq 0$ such that ${\bf z}$ is $n$-regular for all
$n\geq m+1$ is denoted by $a({\bf z})$ (see \cite{trung1},
Section~2). In other words,
\begin{eqnarray*}
a({\bf z})={\rm max}\{ a((z_{1},\ldots ,z_{i-1})F:z_{i}/(z_{1},\ldots ,
z_{i-1})F) \mid i=1,\ldots ,s\}.
\end{eqnarray*}

A sequence ${\bf z}=z_{1},\ldots ,z_{s}$ of homogeneous elements of
$U$ is called a {\em $U_{+}$-filter-regular sequence with respect to
$F$} if $z_{i}\not\in \fp$ for any associated prime ideal $\fp$ of
$F/(z_{1},\ldots ,z_{i-1})F$, $\fp\not\supseteq U_{+}$, for all
$i=1,\ldots ,s$ (see \cite{trung1}, Section~2, \cite{trung2},
Section~2).

\begin{Lemma} {\rm (\cite{bs}, 18.3.8, \cite{trung1}, 2.1)}
Let ${\bf z}=z_{1},\ldots ,z_{s}$ be a sequence of homogeneous
elements of $U$. Then ${\bf z}$ is a $U_{+}$-filter regular sequence
with respect to $F$ if and only if $a({\bf z})<\infty$.
\end{Lemma}

\begin{Lemma}\label{trung1/2.3} {\rm (\cite{trung1}, 2.3)}
Let $z\in U_{1}$ be a homogeneous $U_{+}$-filter-regular element with
respect to $F$. Then, for all $i\geq 0$,
\begin{eqnarray*}
a_{i+1}(F)+1\leq a_{i}(F/zF)\leq \mbox{\rm max}\{ a_{i}(F),
a_{i+1}(F)+1\}.
\end{eqnarray*}
\end{Lemma}

\begin{Lemma}\label{trung2/2.2} {\rm (\cite{trung2}, 2.2)}
Let ${\bf z}=z_{1},\ldots ,z_{s}$ be a $U_{+}$-filter-regular sequence
with respect to $F$, $z_{i}\in U_{1}$ for all $i=1,\ldots ,s$. Then
\begin{eqnarray*}
&&a({\bf z})={\rm max}\{ a_{i}(F)+i\mid i=0,\ldots ,s-1\}\mbox{ and,
for all }0\leq t\leq s,\\ &&{\rm max}\{ a_{i}(F)+i\mid i=0,\ldots
,t\}= {\rm max}\{ a((z_{1},\ldots ,z_{i})F:U_{+}/(z_{1},\ldots
,z_{i})F) \mid i=0,\ldots, t\}.
\end{eqnarray*}
\end{Lemma}

\begin{Proposition}\label{reg-a} {\rm (\cite{trung2}, 2.4)}
Let ${\bf z}=z_{1},\ldots ,z_{s}$ be a $U_{+}$-filter-regular sequence
with respect to $F$, $z_{i}\in U_{1}$, $i=1,\ldots ,s$, which
generates a reduction $Q$ of $U_{+}$ with respect to $F$. Then
\begin{eqnarray*}
\reg{(F)}=\mbox{\rm max}\{ a({\bf z}),{\rm rn}_{Q}(U_{+};F)\}.
\end{eqnarray*}
\end{Proposition}

\demo By Lemma~\ref{trung2/2.2}, $a({\bf z})={\rm max}\{ a((z_{1},\ldots
,z_{i})F:U_{+}/(z_{1},\ldots ,z_{i})F)\mid i=0,\ldots
,s-1\}$. Further, ${\rm rn}_{Q}(U_{+};F)=a(F/QF)=a((z_{1},\ldots
,z_{s})F:U_{+}/(z_{1},\ldots ,z_{s})F)$. Therefore,
\begin{eqnarray*}
&&{\rm max}\{ a({\bf z}),{\rm rn}_{Q}(U_{+};F)\}={\rm max}\{
a((z_{1},\ldots ,z_{i})F:U_{+}/(z_{1},\ldots ,z_{i})F) \mid
i=0,\ldots, s\}=\\&&{\rm max}\{ a_{i}(F)+i\mid i=0,\ldots ,s\}.
\end{eqnarray*}
Since $\reg{(F)}={\rm max}\{ a_{i}(F)+i\mid i\geq 0\}$, it is enough
to show that $H^{i}_{U_{+}}(F)=0$ for all $i>s$. If $s=0$, then $0$ is
a reduction of $U_{+}$ with respect to $F$ and $F_{n}$ for all large
$n$. So $F$ is a $U_{+}$-torsion module and $H^{i}_{U_{+}}(F)=0$ for
all $i>0$ (see e.g. \cite{bs}, 2.1.7). If $s\geq 1$, by induction,
$H_{U_{+}}^{i}(F/z_{1}F)=0$ for all $i>s-1$. So
$a_{i}(F/z_{1}F)=-\infty$ for all $i>s-1$. By Lemma~\ref{trung1/2.3},
$a_{i+1}(F)=-\infty$ and $H^{i+1}_{U_{+}}(F)=0$ for all $i>s$. \qed

\medskip

Now take $A$ a noetherian ring, $I$ an ideal of $A$ and $M$ a finitely
generated $A$-module. Consider $\rees{I}=\oplus _{n\geq
0}I^{n}t^{n}\subset A[t]$ as a subring of $A[t]$.

\begin{Lemma}\label{filter} {\rm (\cite{trung2}, 4.1)}
Let $A$ be a noetherian ring, let $I$ be an ideal of $A$ and let $M$
be a finitely generated $A$-module. Let $x_{1},\ldots ,x_{s}$ be a
sequence of elements in $I$. Then $x_{1}t,\ldots ,x_{s}t$ is a
$\rees{I}_{+}$-filter-regular sequence with respect to $\reesm{I}{M}$
if and only if for all large $n\geq 1$,
\begin{eqnarray*}
[(x_{1},\ldots ,x_{i-1})I^{n}M:x_{i}]\cap I^{n}M=(x_{1},\ldots ,
x_{i-1})I^{n-1}M\mbox{ for }i=1,\ldots ,s. \quad (^{*})
\end{eqnarray*}
If that is the case, $a({\bf z})$ is the least integer $r$ such that
$(^{*})$ holds for all $n\geq r+1$.
\end{Lemma}

\demo ${\bf z}=x_{1}t,\ldots ,x_{s}t$ is a
$\rees{I}_{+}$-filter-regular sequence with respect to $\reesm{I}{M}$
if and only if $[(x_{1}t,\ldots ,x_{i-1}t)\reesm{I}{M}:x_{i}t]_{n}$ is
equal to $[(x_{1}t,\ldots ,x_{i-1}t)\reesm{I}{M}]_{n}$ for all large
$n\geq 1$. But the first module is equal to $[(x_{1},\ldots
,x_{i-1})I^{n}M:x_{i}]\cap I^{n}M$ and the second is equal to
$(x_{1},\ldots , x_{i-1})I^{n-1}M$. \qed

\begin{Proposition}\label{regu} {\rm (\cite{trung2}, 4.2)}
Let $A$ be a noetherian ring, let $I$ be an ideal of $A$ and let $M$
be a finitely generated $A$-module. Let $J=(x_{1},\ldots ,x_{s})$ be a
reduction of $I$ with respect to $M$. Suppose that ${\bf
z}=x_{1}t,\ldots ,x_{s}t$ is a $\rees{I}_{+}$-filter-regular sequence
with repect to $\reesm{I}{M}$. Then
\begin{eqnarray*}
\reg(\reesm{I}{M})= \mbox{\rm min}\{ r\geq 0\mid r\geq {\rm
  rn}_{J}(I;M)\mbox{ and }(^{*})\mbox{ holds for all }n\geq r+1\}.
\end{eqnarray*}
\end{Proposition}

\demo Let $Q=({\bf z})$ denote the ideal generated by ${\bf
z}=x_{1}t,\ldots ,x_{s}t$, $U=\rees{I}$ the Rees algebra of $I$ and
$F=\reesm{I}{M}$ the Rees module of $I$ with respect to $M$. Since $J$
is a reduction of $I$ with respect to $M$, then $Q$ is a reduction of
$U_{+}$ with respect to $F$. Moreover, if $I^{r+1}M=JI^{r}M$, then
$U_{+}^{r+1}F=QU_{+}^{r}F$ and ${\rm rn}_{Q}(U_{+};F)={\rm
rn}_{J}(I;M)$. By Proposition~\ref{reg-a}, $\reg(F)={\rm max}\{ a({\bf
z}),{\rm rn}_{J}(I;M)\}$. The conclusion follows from
Lemma~\ref{filter}. \qed

\section{Relation type and reduction number}\label{rtrn}

The first result of the section suggests the relationship subsisting
between the relation type and the reduction number (see
\cite{vasconcelos3}, page 63).

\begin{Proposition}\label{rtlrtrn}
Let $A$ be a ring, $I$ an ideal of $A$ and $M$ an $A$-module. Let
$J\subset I$ be a reduction of $I$ with respect to $M$ and with
reduction number ${\rm rn}_{J}(I;M)=r$. Then
\begin{eqnarray*}
{\rm rt}(I;M)\leq \rn_{J}(I;M)+\rt(J;I^{r}M).
\end{eqnarray*}
\end{Proposition}

\demo Let us prove that $E(I;M)_{n}=0$ for all $n\geq
r+\rt(J;I^{r}M)+1$.  Write $n=r+k$, where $k\geq \rt(J;I^{r}M)+1(\geq
2)$. In particular, $I^{n}M=J^{k}I^{r}M=JI^{n-1}M$,
$I^{n-1}M=J^{k-1}I^{r}M=JI^{n-2}M$ and
$I^{n-2}M=J^{k-2}I^{r}M=JI^{n-3}M$ (where $I^{n-3}=A$ if $n=2$ and
$r=0$). Consider the following diagram:

\begin{picture}(330,85)(15,0)

\put(100,60){\makebox(0,0){\mbox{\footnotesize
$\Lambda_{2}(J)\otimes J^{k-2}I^{r}M$}}}
\put(200,60){\makebox(0,0){\mbox{\footnotesize
$J\otimes J^{k-1}I^{r}M$}}}
\put(300,60){\makebox(0,0){\mbox{\footnotesize
$J^{k}I^{r}M$}}}
\put(360,60){\makebox(0,0){$0$}}

\put(140,60){\vector(1,0){30}}
\put(155,68){\makebox(0,0){\mbox{\footnotesize $\partial ^{\prime}_{2,k}$}}}
\put(230,60){\vector(1,0){50}}
\put(255,68){\makebox(0,0){\mbox{\footnotesize $\partial ^{\prime}_{1,k}$}}}
\put(320,60){\vector(1,0){30}}

\put(140,20){\vector(1,0){30}}
\put(155,12){\makebox(0,0){\mbox{\footnotesize $\partial _{2,n}$}}}
\put(230,20){\vector(1,0){50}}
\put(255,12){\makebox(0,0){\mbox{\footnotesize $\partial _{1,n}$}}}
\put(320,20){\vector(1,0){30}}

\put(100,20){\makebox(0,0){{\footnotesize $\Lambda_{2}(I)\otimes I^{n-2}M$}}} 
\put(200,20){\makebox(0,0){{\footnotesize $I\otimes I^{n-1}M$}}} 
\put(300,20){\makebox(0,0){{\footnotesize $I^{n}M$}}}
\put(360,20){\makebox(0,0){$0$}} 
\put(365,17){\makebox(0,0){.}}

\put(100,50){\vector(0,-1){20}}
\put(90,40){\makebox(0,0){\mbox{\footnotesize $g$}}}
\put(200,50){\vector(0,-1){20}}
\put(190,40){\makebox(0,0){\mbox{\footnotesize $f$}}}
\put(299,50){\line(0,-1){20}}
\put(301,50){\line(0,-1){20}}
\end{picture}

\noindent The top row represents the last three nonzero terms of the
$k$-th homogeneous part of the Koszul complex induced by the
$\rees{J}$-linear form $J\otimes \reesm{J}{I^{r}M}\rightarrow
\reesm{J}{I^{r}M}$ and the bottom row represents the last three
nonzero terms of the $n$-th homogeneous part of the Koszul complex
induced by the $\rees{I}$-linear form $I\otimes
\reesm{I}{M}\rightarrow \reesm{I}{M}$. The Koszul differentials are
defined as usual (e.g. \cite{bh}, Definition~1.6.1, see also the proof
of Proposition \ref{E(x,y)}). The vertical morphisms are induced by
the inclusion $J\subset I$ and define a morphism of complexes.  By
\cite{planas}, Proposition 2.6, the first homology groups of these
complexes are ${\rm ker}\partial ^{\prime}_{1,k}/{\rm im}\partial
^{\prime}_{2,k}=E(J;I^{r}M)_{k}$ and ${\rm ker}\partial _{1,n}/{\rm
im}\partial _{2,n}=E(I;M)_{n}$. Thus we want to prove ${\rm
ker}\partial _{1,n}\subset {\rm im}\partial _{2,n}$.  Take $u=\sum
_{i}x_{i}\otimes m_{i}\in I\otimes I^{n-1}M$ such that $\partial
_{1,n}(u)=\sum _{i}x_{i}m_{i}=0$. Write each $m_{i}=\sum
_{j}y_{i,j}m_{i,j}$, $y_{i,j}\in J$, $m_{i,j}\in I^{n-2}M$. Take
$v=\sum _{i,j}y_{i,j}\wedge x_{i}\otimes m_{i,j}\in \Lambda
_{2}(I)\otimes I^{n-2}M$. Then $\partial _{2,n}(v)=u-w$, where $w=
\sum _{i,j}y_{i,j}\otimes x_{i}m_{i,j}\in I\otimes I^{n-1}M$.
Consider $w^{\prime}= \sum _{i,j}y_{i,j}\otimes x_{i}m_{i,j}\in
J\otimes J^{k-1}I^{r}M$. Remark that $\partial
^{\prime}_{1,k}(w^{\prime})=\partial _{1,n}(f(w^{\prime}))=\partial
_{1,n}(w)=0$. Since $k\geq {\rm rt}(J;I^{r}M)+1$, then
$E(J;I^{r}M)_{k}=0$ and $w^{\prime}\in {\rm im}\partial
^{\prime}_{2,k}$. Take $t^{\prime}\in \Lambda_{2}(J)\otimes
J^{k-2}I^{r}M$ such that $\partial
_{2,k}(t^{\prime})=w^{\prime}$. Then
$\partial_{2,n}(v+g(t^{\prime}))=u-w+
f(\partial^{\prime}_{2,k}(t^{\prime}))=u-w+f(w^{\prime})=u$ and
$u\in{\rm im}\partial _{2,n}$. \qed

\medskip

The purpose now is to control the relation type of the reduction $J$
with respect to $I^{r}M$. We will use the filter-regular conditions
$(^{*})$ of Lemma~\ref{filter}, which firstly appeared, to our
knowledge, in a paper by Costa for $M=A$ and $J=I$ (\cite{costa}, page
258).

\begin{Proposition}\label{costacond}
Let $A$ be a ring, let $I$ be an ideal of $A$ and let $M$ be an
$A$-module. Let $J=(x_{1},\ldots,x_{s})\subset I$ be a reduction of
$I$ with respect to $M$ and with reduction number ${\rm
rn}_{J}(I;M)=r$. Suppose that there exists $k\geq 1$ such that for all
$n\geq r+k$ and all $i=1,\ldots ,s$,
\begin{eqnarray*}
[(x_{1},\ldots,x_{i-1})I^{n}M:x_{i}]\cap
I^{n}M=(x_{1},\ldots,x_{i-1})I^{n-1}M.
\end{eqnarray*}
Then ${\rm rt}(J;I^{r}M)\leq k$. Moreover, if $A$ is noetherian and
$M$ is finitely generated, then $\rn_{J}(I;M)\leq
\reg(\reesm{I}{M})\leq \rn_{J}(I;M)+k-1$.
\end{Proposition}

\demo Write $J_{0}=0$ and $J_{i}=(x_{1},\ldots,x_{i})$ for
$i=1,\ldots, s$. Let $m\geq k+1$ and consider the last three nonzero
terms of the $m$-th homogeneous part of the Koszul complex induced by
the $\rees{J_{i}}$-linear form $J_{i}\otimes \reesm{J}{I^{r}M}
\rightarrow \reesm{J}{I^{r}M}$:
\begin{eqnarray*}
\Lambda _{2}(J_{i})\otimes J^{m-2}I^{r}M\buildrel
\partial_{2,m-2}\over \longrightarrow J_{i}\otimes
J^{m-1}I^{r}M\buildrel \partial_{1,m-1}\over\longrightarrow
J^{m}I^{r}M\to 0.
\end{eqnarray*}
If $i=s$, then $J_{s}=J$ and one has the Koszul complex
\begin{eqnarray*}
\Lambda _{2}(J)\otimes J^{m-2}I^{r}M\buildrel
\partial_{2,m-2}\over \longrightarrow J\otimes
J^{m-1}I^{r}M\buildrel \partial_{1,m-1}\over\longrightarrow
J^{m}I^{r}M\to 0,
\end{eqnarray*}
whose first homology group ${\rm ker}\partial_{1,m-1}/{\rm
im}\partial_{2,m-2}$ is, by \cite{planas}, Proposition 2.6, equal to
the module of $m$-effective relations $E(J;I^{r}M)_{m}$. Thus, it is
enough to prove by induction on $i=1,\ldots,s$, that ${\rm
ker}\partial _{1,m-1}\subset {\rm im}\partial _{2,m-2}$ for all $m\geq
k+1$ (remark that, in this case, $r+m-1\geq r+k$). If $i=1$, let
$z=x_{1}\otimes c\in J_{1}\otimes J^{m-1}I^{r}M$ such that $0=\partial
_{1,m-1}(z)=x_{1}c$. Then $c\in (0:x_{1})\cap
J^{m-1}I^{r}M=(0:x_{1})\cap I^{r+m-1}M=J_{0}I^{r+m-2}M=0$. Thus
$z=0$. If $i=s$, let $z=\sum_{i=1}^{s}x_{i}\otimes c_{i}\in J\otimes
J^{m-1}I^{r}M$ such that $0=\partial _{1,m-1}(z)=\sum
_{i=1}^{s}x_{i}c_{i}$. Then
$x_{s}c_{s}=-\sum_{i=1}^{s-1}x_{i}c_{i}$. Thus $c_{s}\in
(J_{s-1}J^{m-1}I^{r}M:x_{s})\cap
J^{m-1}I^{r}M=(J_{s-1}I^{r+m-1}M:x_{s})\cap
I^{r+m-1}M=J_{s-1}I^{r+m-2}M$. Thus
$c_{s}=\sum_{i=1}^{s-1}x_{i}\lambda _{i}$, $\lambda _{i}\in
I^{r+m-2}M=J^{m-2}I^{r}M$. Take $u=\sum_{i=1}^{s-1}x_{i}\otimes
(c_{i}+x_{s}\lambda_{i})\in J_{s-1}\otimes J^{m-1}I^{r}M$. Then
$\partial _{1,m-1}(u)=0$. By induction hypothesis, there exists $v\in
\Lambda _{2}(J_{s-1})\otimes J^{m-2}I^{r}M$ such that $\partial
_{2,m-2}(v)=u$. Take $w=v+\sum_{i=1}^{s-1}(x_{j}\wedge x_{s})\otimes
\lambda _{i}$ and one has $\partial _{2,m-2}(w)=z$. This proves ${\rm
rt}(J;I^{r}M)\leq k$. The second assertion follows from
Proposition~\ref{regu}. \qed

\medskip

Let $A$ be a noetherian ring, $J\subset I$ two ideals of $A$ and $M$ a
finitely generated $A$-module. Let $x_{1},\ldots ,x_{s}$ be a system
of generators of $J$. Write, as before, $J_{0}=0$ and
$J_{i}=(x_{1},\ldots ,x_{i})$ for $i=1,\ldots ,s$. The sequence
$x_{1},\ldots ,x_{s}$ is said to be a {\em $d$-sequence with respect
to $M$} if any $x_{j}$ is not contained in the ideal generated by the
others $x_{i}$ and for all $k\geq i+1$ and all $i\geq 0$,
$(J_{i}M:x_{i+1}x_{k})=(J_{i}M:x_{k})$. It is known that this last
condition is equivalent to $(J_{i}M:x_{i+1})\cap JM=J_{i}M$ for all
$i=0,\ldots ,s-1$. Let $\agrm{I}{M}$ be the associated graded module
of $I$ with respect to $M$ and denote by $x_{1}^{*},\ldots ,x_{s}^{*}$
the images of $x_{1},\ldots,x_{s}$ in $I/I^{2}\subset \agr{I}$. The
sequence $x_{1},\ldots ,x_{s}$ is said to be a {\em complete
$d$-sequence with respect to $I$ and $M$} if $x_{1},\ldots ,x_{s}$ is
a $d$-sequence with respect to $M$ and $x_{1}^{*},\ldots ,x_{s-1}^{*}$
is a $\agrm{I}{M}$-regular sequence (see \cite{huckaba3} and
\cite{trung2}). If $A$ is local, it can be shown that
$x_{1}^{*},\ldots ,x_{s-1}^{*}$ is a $\agrm{I}{M}$-regular sequence if
and only if $x_{1},\ldots ,x_{s-1}$ is an $M$-regular sequence and,
for all $n\geq 0$ and all $i=1,\ldots ,s-1$, the $n$-th
Valabrega-Valla module $VV_{J_{i}}(I;M)_{n}=J_{i}M\cap
I^{n+1}M/J_{i}I^{n}M$ is equal to zero (see e.g. \cite{huckaba2},
Lemma~2.2, \cite{cz}, Proposition~2.3).

Huckaba proved that if $A$ is noetherian local, if $I$ is an ideal
with analytic spread $l(I)$ equal to the height of the ideal ${\rm
ht}(I)$ or ${\rm ht}(I)+1$ and with a minimal reduction $J$ generated
by a complete $d$-sequence with respect to $I$, then ${\rm rt}(I)\leq
{\rm rn}_{J}(I)+1$ (see \cite{huckaba2}, Theorem~2.3 and
\cite{huckaba3}, Theorem~1.4). Later, Trung proved that, in general,
$\rt(I)\leq \reg(\rees{I})+1$ and that if $I$ has a a reduction $J$
generated by a complete $d$-sequence with respect to $I$, then
$\reg(\rees{I})=\rn_{J}(I)$ (see \cite{trung2}, Proposition~2.6 and
Theorem~6.4; for more related results on this topic see also
\cite{schenzel} and \cite{trung1}).  From our
Propositions~\ref{rtlrtrn} and \ref{costacond}, we obtain a
generalization of these results. Our proof closely follows ideas of
Trung in \cite{trung2}.

\begin{Theorem}\label{cds}
Let $A$ be a noetherian ring, let $I$ be an ideal of $A$ and let $M$
be a finitely generated $A$-module. Let $J=(x_{1},\ldots
,x_{s})\subset I$ be a reduction of $I$ with respect to $M$ and with
reduction number ${\rm rn}_{J}(I;M)=r$. Suppose that
\begin{itemize}
\item[$(i)$] $x_{1},\ldots,x_{s}$ is a $d$-sequence with respect to $M$.
\item[$(ii)$] $x_{1},\ldots,x_{s-1}$ is an $M$-regular sequence.
\item[$(iii)$] $(x_{1},\ldots ,x_{i})M\cap I^{r+1}M=(x_{1},\ldots ,
x_{i})I^{r}M$ for all $i=1,\ldots,s-1$.
\end{itemize}
Then $\rt(J;I^{r}M)=1$, $\rt(I;M)\leq \rn_{J}(I;M)+1$ and
$\rn_{J}(I;M)=\reg(\reesm{I}{M})$.
\end{Theorem}

\demo Write $J_{0}=0$ and $J_{i}=(x_{1},\ldots,x_{i})$ for
$i=1,\ldots, s$. Using $(ii)$ and $(iii)$, we obtain
$(J_{i-1}M:x_{i})\cap I^{r+1}M=J_{i-1}M\cap I^{r+1}M= J_{i-1}I^{r}M$
for $i=1,\ldots ,s-1$. By $(i)$, $(J_{s-1}M:x_{s})\cap
JM=J_{s-1}M$. Since $I^{r+1}M=JI^{r}M$, then $(J_{s-1}M:x_{s})\cap
I^{r+1}M=J_{s-1}M\cap I^{r+1}M$ which, by $(iii)$, is equal to
$J_{s-1}I^{r}M$. Thus, for all $i=1,\ldots ,s$,
\begin{eqnarray*}
(J_{i-1}M:x_{i})\cap I^{r+1}M=J_{i-1}I^{r}M.
\end{eqnarray*}
A straighforward generalization to modules of Theorem 4.8, $(i)$ in
\cite{trung2}, allows us to assert that for all integers $n\geq r+1$
and for all $i=1,\ldots ,s$,
\begin{eqnarray*}
(J_{i-1}M:x_{i})\cap I^{n}M=J_{i-1}I^{n-1}M,
\end{eqnarray*}
which clearly implies for all integers $n\geq r+1$ and for all
$i=1,\ldots ,s$,
\begin{eqnarray*}
(J_{i-1}I^{n}M:x_{i})\cap I^{n}M=J_{i-1}I^{n-1}M.
\end{eqnarray*}
By Proposition~\ref{rtlrtrn}, $\rt(I;M)\leq \rn_{J}(I;M)+{\rm
rt}(J;I^{r}M)$ and, by Proposition~\ref{costacond}, $\rt(I;J^{r}M)=1$
and $\rn_{J}(I;M)=\reg(\reesm{I}{M})$.  \qed

\section{Integral degree of a ring}\label{id}

In this section we introduce the {\em integral degree}, an invariant
associated to the ring, which later will be used to bound the
reduction number. Let $A\subset B$ be a ring extension. Recall that an
element $b\in B$ is said to be integral over $A$ if there exist
$a_{i}\in A$ and an integral equation of degree $n\geq 1$:
\begin{eqnarray*}
b^{n}+a_{1}b^{n-1}+a_{2}b^{n-2}+\ldots +a_{n-1}b+a_{n}=0.
\end{eqnarray*}
If $b\in B$ is integral over $A$, we will call the {\em integral
degree of} $b$ {\em over} $A$ to the integer:
\begin{eqnarray*}
{\rm id}_{A}(b)={\rm min} \{ n\geq 1\mid b\mbox{ satisfies an integral
equation of degree }n\}.
\end{eqnarray*}
Let $A\subset C\subset B$, $C$ an $A$-submodule of $B$. Suppose the
elements of $C$ are integral over $A$.  Then the {\em integral
degree of} $C$ {\em over} $A$ is defined as the integer (possibly
infinite):
\begin{eqnarray*}
\ud_{A}(C)={\rm sup}\, \{ {\rm id}_{A}(c)\mid c\in C\}.
\end{eqnarray*}
Remark that $\ud_{A}(C)=1$ if and only if $A=C$.

As usual, $\mu _{A}(\cdot )$ stands for the minimal number of
generators as an $A$-module.

\begin{Proposition}\label{genatiyah}
Let $A\subset B$ be a ring extension, $b\in B$ and $n\geq 1$. Then the
following conditions are equivalent:
\begin{itemize}
\item[$(i)$] $b$ is integral over $A$ and ${\rm id}_{A}(b)\leq n$.
\item[$(ii)$] $A[b]$ is a finitely generated $A$-module and $\mu
_{A}(A[b])\leq n$.
\item[$(iii)$] There exists a ring $C$, $A\subset A[b]\subset C\subset
B$, such that $C$ is a finitely generated $A$-module and $\mu
_{A}(C)\leq n$.
\item[$(iv)$] There exists a faithful $A[b]$-module $M$ such that $M$
is a finitely generated $A$-module and $\mu _{A}(M)\leq n$.
\end{itemize}
\end{Proposition}

\demo It follows from \cite{am}, Proposition~5.1, just taking into
account the definition of ${\rm id}_{A}(b)$. \qed

\medskip

\noindent Next we prove that the integral degree of the sum or product
of two integral elements is, in fact, bounded above by the product of
their integral degrees.

\begin{Corollary}\label{nelements}
Let $A\subset B$ be a ring extension and $b_{1},\ldots ,b_{n}\in B$
integral over $A$. Then $A[b_{1},\ldots ,b_{n}]$ is a finitely
generated $A$-module, $A\subset A[b_{1},\ldots ,b_{n}]$ is an integral
extension and:
\begin{eqnarray*}
{\rm max}\{ {\rm id}_{A}(b_{i})\}\leq \ud_{A}(A[b_{1},\ldots
,b_{n}])\leq \mu _{A}(A[b_{1},\ldots ,b_{n}])\leq \prod_{i=1}^{n}{\rm
id}_{A}(b_{i}).
\end{eqnarray*}
In particular, if $b\in B$ is integral over $A$, then ${\rm
id}_{A}(b)=\ud_{A}(A[b])=\mu _{A}(A[b])$.
\end{Corollary}

\demo Let $C=A[b_{1},\ldots ,b_{n}]$ and $m=\prod_{i=1}^{n}{\rm
id}_{A}(b_{i})$. Then it is clear that $\mu _{A}(C)\leq m$. Now take
any $b\in C\subset B$. So we have $A\subset A[b]\subset C\subset B$
with $\mu _{A}(C)=r\leq m$. By Proposition \ref{genatiyah},
$(iii)\Rightarrow (i)$, $b$ is integral over $A$ and ${\rm
id}_{A}(b)\leq r$ and taking the supremum over all $b\in C$, then
$\ud_{A}(C)\leq r=\mu_{A}(C)$. As for the second assertion, just
take $n=1$. \qed

\begin{Corollary}\label{uidlmu}
Let $A\subset B$ be a ring extension. If $B$ is a finitely generated
$A$-module, then $A\subset B$ is integral and
\begin{eqnarray*}
\ud_{A}(B)\leq\mu _{A}(B).
\end{eqnarray*}
\end{Corollary}

\demo If $b\in B$, take $A\subset A[b]\subset B$, with $B$ a finitely
generated $A$-module. By Proposition~\ref{genatiyah},
$(iii)\Rightarrow (i)$, $b$ is integral over $A$ and $\id_{A}(b)\leq
\mu _{A}(B)$. Taking the supremum, $\ud_{A}(B)\leq \mu _{A}(B)$. \qed

\medskip

Let $A$ be noetherian domain and let $\overline{A}$ be the integral
closure of $A$ in its quotient field. If ${\rm dim}\, A\leq 2$, then
$\overline{A}$ is noetherian (see e.g. \cite{matsumura}, 11.7 and
\cite{nagata}, 33.12). Nevertheless, $\overline{A}$ may be a non
finitely generated $A$-module, as an example of Akizuki shows
(\cite{akizuki} or \cite{reid}, 9.5). Next, we want to prove that the
ring $A$ in the example of Akizuki has at least finite integral degree
$\ud_{A}(\overline{A})$. Before that, and for easy reference, we state
the following lemma.

\begin{Lemma}\label{idrtrn}
Let $A$ be a ring and $x,y\in A$, with $x$ regular. The following are
equivalent.
\begin{itemize}
\item[$(i)$] $y/x$ is integral over $A$ and $\id_{A}(y/x)\leq n$.
\item[$(ii)$] $(x)$ is a reduction of $(x,y)$ and ${\rm
  rn}_{(x)}(x,y)\leq n-1$.
\item[$(iii)$] $x(x,y)^{n-1}:y^{n}=A$.
\end{itemize}
In particular, if $y/x$ is integral over $A$, then $\id_{A}(y/x)={\rm
rt}(x,y)={\rm rn}_{(x)}(x,y)+1$
\end{Lemma}

\proof Take $y/x\in\overline{A}$ with $\id_{A}(y/x)\leq n$. Then,
there exist $a_{i}\in A$ such that $(y/x)^{n}+a_{1}(y/x)^{n-1}+\ldots
+a_{n}=0$. Multiplying by $x^{n}$, one has $y^{n}\in xI^{n-1}$, where
$I=(x,y)$. Thus $I^{n}=xI^{n-1}$, $J=(x)$ is a reduction of $I$ and
${\rm rn}_{J}(I)\leq n-1$. If $J=(x)$ is a reduction of $I=(x,y)$ with
${\rm rn}_{J}(I)\leq n-1$, then $I^{n}=xI^{n-1}$ and $y^{n}\in
xI^{n-1}$. Thus $1\in xI^{n-1}:y^{n}$ and $xI^{n-1}:y^{n}=A$. Finally,
if $xI^{n-1}:y^{n}=A$, where $I=(x,y)$, then $y^{n}\in xI^{n-1}$ and
$y^{n}=b_{1}xy^{n-1}+\ldots +b_{n}x^{n}$, for some $b_{i}\in
A$. Dividing by $x^{n}$ one obtains an integral equation of $y/x$ over
$A$ of degree $n$. In particular, if $y/x\in \overline{A}$ with ${\rm
id}_{A}(y/x)=n\geq 2$, then $xI^{n-2}:y^{n-1}\varsubsetneq
A=xI^{n-1}:y^{n}$, where $I=(x,y)$. By Proposition~\ref{E(x,y)},
$E(I)_{n}\neq 0$ and $E(I)_{n+s}=0$ for all $s\geq 1$. Thus ${\rm
rt}(I)=n={\rm id}_{A}(y/x)$. Moreover, $(i)\Leftrightarrow (ii)$ says
that $J=(x)$ is a reduction of $I=(x,y)$ and that ${\rm
rn}_{J}(I)=n-1$. \qed

\medskip

Now, let us prove that the example of Akizuki has finite integral
degree. Denote by $e(A)$ the multiplicity of $A$.

\begin{Proposition}\label{akizuki} Let $A$ be a one-dimensional 
noetherian local ring. Then 
\begin{eqnarray*}
\ud_{A}(\overline{A})\leq e(A/H^{0}_{\mathfrak{m}}(A))+{\rm
length}(H^{0}_{\mathfrak{m}}(A)).
\end{eqnarray*}
In particular, if $A$ is a domain (as Akizuki's example is), then
$\ud_{A}(\overline{A})\leq e(A)$.
\end{Proposition}

\demo Take $y/x\in \overline{A}$ with $x,y\in A$, $x$ regular, and
$I=(x,y)$ the ideal of $A$ generated by $x,y$. By Lemma~\ref{idrtrn},
$\id_{A}(y/x)={\rm rt}(I)$.  By \cite{planas}, Lemma~6.1, ${\rm
rt}(I)\leq {\rm rt}(IA/J)+{\rm length}(J)$, where
$J=H^{0}_{\mathfrak{m}}(A)$. By \cite{planas}, Lemma~6.3, ${\rm
rt}(IA/J)\leq e(A/J)$. \qed

\medskip

We next see that there exist one-dimensional noetherian domains with
infinite integral degree. Remark that the ring in this example must be
not local nor excellent so that one can not apply the existence of a
uniform bound for the relation type of all ideals (see \cite{planas},
Proposition~6.5 and Theorem~3). The next example is due to Sally and
Vasconcelos (see \cite{sv}, Example 1.4, and also \cite{planas},
Remark 7.3).

\begin{Example}{\rm 
There exists one-dimensional noetherian domains $A$ with
$\ud_{A}(\overline{A})$ infinite.}\end{Example}

\demo Let $t_{1},t_{2},t_{3},\ldots $ be infinitely many
indeterminates over a field $k$. Let $R$ be defined as
$R=k[t_{1}^{2},t_{1}^{3},t_{2}^{3},t_{2}^{4}, t_{2}^{5},\ldots
,t_{n}^{n+1},t_{n}^{n+2},\ldots ,t_{n}^{2n+1},\ldots ]$. Take
$\fp_{n}=(t_{n}^{n+1},t_{n}^{n+2},\ldots ,t_{n}^{2n+1})$, which is a
prime ideal of height 1. Let $S$ be the multiplicative closed set
$R-\cup \fp_{n}$ and $A=S^{-1}R$. One can prove that $A$ is a
one-dimensional noetherian domain and that $t_{n}^{n+2}/t_{n}^{n+1}$
is in $\overline{A}$ and has integral degree $n$.  Therefore
$\ud_{A}(\overline{A})=\infty$. \qed

\medskip

We now give two more properties of the integral degree.

\begin{Proposition}
Let $A\subset B$ and $B\subset C$ be integral extensions. Then
$A\subset C$ is an integral extension and
\begin{eqnarray*}
\ud_{A}(C)\leq \ud_{A}(B)^{{\rm d}_{B}(C)}\cdot \ud_{B}(C).
\end{eqnarray*}
\end{Proposition}

\demo If $c\in C$, there exists an equation $c^{n}+b_{1}c^{n-1}+\ldots
+b_{n-1}c+b_{n}=0$, with $b_{i}\in B$ and $n\leq
\ud_{B}(C)$. Take $D=A[b_{1},\ldots ,b_{n}]$. Since $A\subset B$ is an
integral extension, all $b_{i}$ are integral over $A$ and, by
Corollary \ref{nelements}, $D$ is a finitely generated $A$-module and
$\mu _{A}(D)\leq \prod_{i=1}^{n}{\rm id}_{A}(b_{i})\leq
\ud_{A}(B)^{{\rm d}_{B}(C)}$. On the other hand, $c$ is clearly
integral over $D$ and $D[c]$ is a finitely generated $D$-module with
$\mu _{D}(D[c])\leq n\leq \ud_{B}(C)$. Since $D$ is a finitely
generated $A$-module and $D[c]$ is a finitely generated $D$-module,
then $D[c]$ is a finitely generated $A$-module. So we have $A\subset
A[c]\subset D[c]\subset C$ with $D[c]$ a finitely generated $A$-module
with $\mu _{A}(D[c])\leq \mu _{A}(D)\mu _{D}(D[c])\leq
\ud_{A}(B)^{{\rm d}_{B}(C)}\cdot \ud_{B}(C)$. Applying Proposition
\ref{genatiyah}, $(iii)\Rightarrow (i)$, we deduce that $c$ is
integral over $A$ and $\id_{A}(c)\leq \ud_{A}(B)^{{\rm d}_{B}(C)}\cdot
\ud_{B}(C)$. \qed

\begin{Proposition}\label{uidlocalitzat}
Let $A\subset B$ be an integral extension and $S$ a multiplicatively
closed subset of $A$. Then $S^{-1}A\subset S^{-1}B$ is an integral
extension and
\begin{eqnarray*}
\ud_{S^{-1}A}(S^{-1}B)\leq \ud_{A}(B).
\end{eqnarray*}
In particular, if $A$ is reduced, $\mathfrak{p}$ is a prime ideal of
$A$ and $\overline{A}$ and $\overline{A_{\fp}}$ are the integral
closures of $A$ and $A_{\fp}$ in their total quotient rings, then
\begin{eqnarray*}
\ud_{A_{\fp}}(\overline{A_{\fp}})\leq \ud_{A}(\overline{A}).
\end{eqnarray*}
\end{Proposition}

\demo Let $b/s\in S^{-1}B$, $b\in B$, $s\in S$. Then
$b^{n}+a_{1}b^{n-1}+\ldots +a_{n-1}b+a_{n}=0$, say, and multiplying by
$s^{-n}$ we get $(b/s)^{n}+(a_{1}/s)(b/s)^{n-1}+\ldots
+(a_{n-1}/s^{n-1})(b/s)+a_{n}/s^{n}=0$. So $b/s$ is integral over
$S^{-1}A$ and ${\rm id}_{S^{-1}A}(b/s)\leq {\rm id}_{A}(b)\leq
\ud_{A}(B)$. If $A$ is reduced, then
$S^{-1}\overline{A}=\overline{S^{-1}A}$ (see e.g. \cite{hw}, Lemma
2.1). Therefore, if $S=A-\fp$,
$\overline{(A_{\fp})}=\overline{S^{-1}A}=S^{-1}\overline{A}$ and
$\ud_{A_{\fp}}(\overline{A_{\fp}})=\ud_{S^{-1}A}(S^{-1}\overline{A})\leq
\ud_{A}(\overline{A})$. \qed

\medskip

If $A$ is not reduced, Proposition \ref{uidlocalitzat} may fail. The
next example is taken from \cite{hw}.

\begin{Example} {\rm 
Let $k$ be a field and $A=k\lbrack\!\lbrack x,y,z\rbrack\!\rbrack
/(x^{3}-y^{2})(x,y,z)$. Since the maximal ideal annihilates the
non-zero element $x^{3}-y^{2}$, $A$ is integrally closed and
$\ud_{A}(\overline{A})=1$. Set $S$ the multiplicatively closed set
$\{z^{n},n\geq 0\}$. Then $S^{-1}A= K\lbrack\!\lbrack
x,y\rbrack\!\rbrack /(x^{3}-y^{2})$, where $K$ is the quotient field
of $k\lbrack\!\lbrack z\rbrack\!\rbrack$, and one can prove that
$\ud_{S^{-1}A}(\overline{S^{-1}A})=2$.}
\end{Example}

\section{Integral degree and reduction number}\label{idrn}

We now prove the main result of the paper. Recall that if $I$ is a
regular ideal having principal reductions $J_{1}$ and $J_{2}$ with
${\rm rn}_{J_{1}}(I)=n$ and ${\rm rn}_{J_{2}}(I)=m$, then Huckaba
proved that $n=m$ (see \cite{huckaba1}, where the local assumption is
not needed; it could also be deduced from Theorem~\ref{cds}). We will
denote ${\rm rn}(I)$ to the $J$-reduction number of $I$ for any
principal reduction $J$ of $I$.

\begin{Theorem}\label{id=rn}
Let $A$ be a noetherian ring, $A\supset\bq$. Then
\begin{eqnarray*}
\ud_{A}(\overline{A})={\rm sup}\, \{ {\rm rn}(I)\mid I\mbox{ regular
ideal of }A\mbox{ having a principal reduction}\} +1.
\end{eqnarray*}
\end{Theorem}

\demo Set $\sigma={\rm sup}\, \{\, \rn(I)\mid I \mbox{ regular ideal
of }A\mbox{ having a principal reduction}\}+1$ and
$d=\ud_{A}(\overline{A})$. Take $I$ any regular ideal of $A$ having a
principal reduction $J=(x)$, which is also regular. Then
$I^{n+1}=xI^{n}$ for some $n\geq 0$.  Set $H=x^{-1}I$. Then $H$ is a
fractional ideal of $A$ with $H^{n+1}=H^{n}$. If $y\in I$,
$(y/x)H^{n}\subset H^{n+1}=H^{n}$. Thus $H^{n}$ is a faithful
$A[y/x]$-module. By Proposition \ref{genatiyah}, $y/x$ is integral
over $A$. Thus ${\rm id}_{A}(y/x)\leq d$. By Lemma~\ref{idrtrn},
$x(x,y)^{d-1}:y^{d}=A$ and $y^{d}\in x(x,y)^{d-1}\subseteq xI^{d-1}$.
Therefore $I^{[d]}\subset xI^{d-1}$, where $I^{[d]}$ stands for the
ideal generated by the $d$-th powers of all elements of $I$.  If
$A\supset \bq$, then $I^{[d]}=I^{d}$ (see e.g. \cite{bourbaki}, A1,
\S~8, n$^{\circ}$~2, page~95). Thus $\rn(I)\leq d-1$ and $\sigma\leq
d$. Now take $x,y\in A$, with $x$ regular, such that $y/x$ is integral
over $A$. By Lemma~\ref{idrtrn}, $\id_{A}(y/x)=\rn(x,y)+1\leq
\sigma$. Therefore $d\leq \sigma$. \qed

\begin{Remark}{\rm Let $A$ be a noetherian ring, $A\supset\bq$.
If $I$ is a regular ideal of $A$ having a principal reduction, by
Theorem~\ref{cds}, $\rt(I)\leq\rn(I)+1$ and
$\reg(\rees{I})=\rn(I)$. Moreover, by Lemma~\ref{idrtrn},
$\id_{A}(y/x)=\rt(x,y)$ for any $x,y\in A$, with $x$ regular and such
that $y/x$ is integral over $A$. In other words,
$\ud_{A}(\overline{A})$ is less than or equal to the supremum of the
relation type of two-generated regular ideals of $A$ having principal
reductions. Therefore, in Theorem~\ref{id=rn}, one can replace
$\rn(I)+1$ by else $\reg(\rees{I})+1$ or else $\rt(I)$. In addition,
one can take the supremum just over the set of two-generated regular
ideals having principal reductions. }\end{Remark}

We state a particular version of Theorem~\ref{id=rn} which will be
used later. Note that here we do not need the hypothesis
$A\supset\bq$.

\begin{Proposition}\label{rtm}
Let $(A,\fm)$ be a noetherian local ring with infinite residue
field. Then
\begin{eqnarray*}
\ud_{A}(\overline{A})={\rm sup}\, \{ {\rm rt}_{\mathfrak{m}}(I)\mid
I\mbox{ two-generated regular ideal of $A$}\}.
\end{eqnarray*}
\end{Proposition}

\demo Set $\sigma={\rm sup}\, \{\, \rt_{\mathfrak{m}}(I)\mid I \mbox{
two-generated regular ideal of }A \}$ and $d=\ud_{A}(\overline{A})$.
Take $I$ a two-generated regular ideal of $A$. Since $A$ is noetherian
local with infinite residue field, $I$ has a minimal reduction $J$
generated by as many elements as the analytic spread $l(I)$ of $I$
(see \cite{nr} or \cite{hs}). If $l(I)=1$, by Remark~\ref{coef-mod}
and Theorems~\ref{cds} and \ref{id=rn}, $\rt_{\mathfrak{m}}(I)\leq
\rt(I)\leq \rn(I)+1\leq d$. If $l(I)=2$, then $I$ is generated by two
analytically independent elements $x,y$ and the fiber cone of $I$,
$\reesq{\fm}{I}=\fcone{I}=\oplus _{n\geq 0}I^{n}/\fm I^{n}$ is
isomorphic to a polynomial ring $(A/\fm)[X,Y]$. Thus ${\rm
rt}_{\mathfrak{m}}(I)={\rm rt}(\fcone{I})=1\leq d$ and $\sigma\leq
d$. Now take $x,y\in A$, with $x$ regular, such that $y/x$ is integral
over $A$. Set $\id_{A}(y/x)=n$. By Lemma~\ref{idrtrn},
$xI^{n-2}:y^{n-1}\varsubsetneq xI^{n-1}:y^{n}=A$. By
Proposition~\ref{E(x,y)}, $E_{\fm}(I)_{n}=A/\fm$ and
$E_{\fm}(I)_{n+s}=0$ for all $s\geq 1$. Thus ${\rm
rt}_{\fm}(I)=n$. Therefore ${\rm id}_{A}(y/x)=n={\rm rt}_{\fm}(I)\leq
\sigma$. Thus $d\leq \sigma$. \qed

\medskip

\begin{Remark}\label{necprin} {\rm 
Clearly, Theorem~\ref{id=rn} is no longer true for ideals having
reductions generated by regular sequences of length $l\geq 2$. For
instance, in the power series ring $A=k\lbrack\!\lbrack
x,y\rbrack\!\rbrack$ over a field $k$, the ideals $I_{n}=(x^{n},y^{n},
x^{n-1}y)$ have reductions $(x^{n},y^{n})$ with reduction number $n-1$
(see \cite{huneke}, Remark 5.8). However, $(x^{n},y^{n})$ does not
verify condition $(iii)$ of Theorem~\ref{cds}.  }\end{Remark}

\section{Uniform Artin-Rees numbers}\label{uarp}

We now can prove all the results related to Artin-Rees properties.

\begin{Theorem}\label{modules}
Let $A$ be a noetherian ring with finite integral degree $\ud
_{A}(\overline{A})=d$. Suppose that $A\supset\bq$. Let $N\subset M$ be
two finitely generated $A$-modules. Let $I$ be a regular ideal of $A$
having a principal reduction generated by a $d$-sequence with respect
to $M/N$. Then, for every integer $n\geq d$,
\begin{eqnarray*}
I^{n}M\cap N=I^{n-d}(I^{d}M\cap N).
\end{eqnarray*}
\end{Theorem}

\demo Since $I$ is regular and has a principal reduction, by
Theorem~\ref{id=rn}, $\rn(I)\leq d-1$. It is enough to prove that
$s(N,M;I)\leq \rn(I)+1$. Let $J=(x)$ be a principal reduction of $I$,
set $r=\rn(I)$ and take $k\geq 1$. Then $I^{r+k}M\cap
N=x^{k}I^{r}M\cap N$. Since $x$ is a $d$-sequence with respect to
$M/N$, $x^{k}I^{r}M\cap N=x^{k-1}(xI^{r}M\cap N)\subseteq
I^{k-1}(I^{r+1}M\cap N)$. Thus $s(N,M;I)\leq r+1$. \qed

\medskip

\begin{Remark}{\rm 
By Proposition \ref{slrt}, $s(N,M;I)\leq {\rm rt}(I;M/N)$. If $I$ has
a principal reduction $J$ generated by a $d$-sequence with respect to
$M/N$, then $J$ is also a principal reduction of $I$ with respect to
$M/N$ and, by Theorem~\ref{cds}, ${\rm rt}(I;M/N)\leq {\rm
rn}_{J}(I;M/N)+1\leq {\rm rn}_{J}(I)+1$. Since $I$ is a regular ideal
having a principal reduction, by Theorem~\ref{id=rn}, ${\rm rn}(I)\leq
d-1$. Therefore, $s(N,M;I)\leq d$, which also proves
Theorem~\ref{modules}.  }\end{Remark}

Our ideal-theoretic version of Theorem~\ref{modules} is the following.

\begin{Theorem}\label{ideals}
Let $A$ be a noetherian ring, $A\supset\bq$. Let $\fa$ be an ideal of
$A$ such that $A/\fa$ has finite integral degree
$\ud_{A/\fa}(\overline{A/\fa})=d$. Let $I$ be an ideal of $A$ such
that $IA/\fa$ has an $A/\fa$-regular principal reduction.
Then, for every integer $n\geq d$,
\begin{eqnarray*}
I^{n}\cap \fa=I^{n-d}(I^{d}\cap \fa).
\end{eqnarray*}
\end{Theorem}

\demo By Proposition~\ref{slrt}, $s(\fa,A;I)\leq {\rm
rt}(I;A/\fa)$. By Remark~\ref{coef-mod}, ${\rm rt}(I;A/\fa)={\rm
rt}(IA/\fa)={\rm rt}(IA/\fa;A/\fa)$. Since $IA/\fa$ is $A/\fa$-regular
and has principal reduction $JA/\fa$, by Theorem~\ref{cds}, ${\rm
rt}(IA/\fa;A/\fa)\leq {\rm rn}_{JA/\fa}(IA/\fa;A/\fa)+1={\rm
rn}_{JA/\fa}(IA/\fa)+1$. By Theorem~\ref{id=rn}, ${\rm
rn}_{JA/\fa}(IA/\fa)\leq d-1$. So $s(\fa,A;I)\leq d$. \qed

\medskip

As a corollary of Theorem~\ref{ideals} we obtain a particular version
of the main result in \cite{planas}. 

\begin{Corollary}\label{cor1dim}
Let $(A,\fm)$ be a noetherian local ring, $A\supset\bq$. Let $\fa$ be
an ideal of $A$ such that $A/\fa$ has finite integral degree
$\ud_{A/\fa}(\overline{A/\fa})=d$. Suppose that ${\rm dim}(A/\fa)\leq
1$. Then, for every integer $n\geq d$ and for every ideal $I$ of $A$
such that $IA/\fa$ is $A/\fa$-regular,
\begin{eqnarray*}
I^{n}\cap \fa=I^{n-d}(I^{d}\cap \fa).
\end{eqnarray*}
\end{Corollary}

\demo Since ${\rm dim}(A/\fa)\leq 1$, every ideal $I$ of $A$ is such
that $IA/\fa$ has a principal reduction. Then apply
Theorem~\ref{ideals}. \qed

\medskip

Remark that by a result of Krull, if $(R,\fn)$ is a noetherian local
non-reduced ring such that $\fn$ contains a regular element, then the
integral closure $\overline{R}$ is not a finite $R$-module (see
e.g. \cite{matsumura},~\S33). In particular, in Theorem~\ref{ideals}
and in Corollary~\ref{cor1dim} (as well as in Theorem~\ref{ideals2}),
setting $R=A/\fa$, if $A/\fa$ has a finite integral closure and
$\fm/\fa$ has a regular element, one deduces that $\fa$ is forced to
be a radical ideal.

The next example, taken from Eisenbud and Hochster in \cite{eh}, shows
that if the integral degree is not finite, then the conclusion of
Theorem~\ref{ideals} may be false.

\begin{Example}\label{ehe}{\rm 
There exist $A$, a two-dimensional noetherian domain, $\fa$, a prime
ideal of $A$, and $\{I_{n}\}_{n}$, a family of two-generated ideals of
$A$ such that $I_{n}A/\fa$ has an $A/\fa$-regular principal reduction,
but, for every integer $n\geq 1$,
\begin{eqnarray*}
I_{n}^{n}\cap \fa\varsupsetneq I_{n}(I_{n}^{n-1}\cap \fa) .
\end{eqnarray*}}
\end{Example}

\demo Let $k$ be an algebraically closed field and $\{ X_{n}\}$, $\{
Y_{n}\}$ two countable families of indeterminates. Set
$f_{n}=X^{n}_{n}-Y^{n+1}_{n}$ and $I_{n}$ the ideal in
$T_{n}=k[X_{1},Y_{1},\ldots ,X_{n},Y_{n}]$ generated by
$f_{2}-f_{1},\ldots ,f_{n}-f_{1}$. Set $S_{n}=T_{n}/I_{n}$ and
$U_{n}=S_{n}-\cup _{i=1}^{n}(X_{i},Y_{i})S_{n}$. $U_{n}$ is a
multiplicatively closed subset of $S_{n}$. Set
$A_{n}=U_{n}^{-1}S_{n}$, $A=\botrel{\lim}{\rightarrow}A_{n}$ and
$x_{n},y_{n}$ and $f$ the images of $X_{n},Y_{n}$ and $f_{n}$ in
$A$. Then $A$ is a two-dimensional noetherian regular factorial ring
whose maximal ideals $I_{n}=(x_{n},y_{n})$ form a countable set. Their
intersection $\cap _{n}I_{n}$ is a prime principal ideal $\fa=(f)$
whose generator $f$ is in $I_{n}^{n}$.  Then $I_{n}A/\fa$ is
$A/\fa$-regular and $y_{n}A/\fa$ is a principal reduction of
$I_{n}A/\fa$. Moreover,
\begin{eqnarray*}
I_{n}(I_{n}^{n-1}\cap \fa)=I_{n}\fa\varsubsetneq \fa =I_{n}^{n}\cap \fa.
\end{eqnarray*}
In particular, by Theorem~\ref{ideals},
$\ud_{A/\fa}(\overline{A/\fa})=\infty$. \qed

\medskip

We now prove that there exists a uniform Artin-Rees modulo $\fm$
number for the set of two-generated regular ideals. Here, we do not
need $A\supset\bq$.

\begin{Theorem}\label{ideals2}
Let $(A,\fm)$ be a noetherian local ring with infinite residue
field. Let $\fa$ be an ideal of $A$ such that $A/\fa$ has finite
integral degree $\ud_{A/\fa}(\overline{A/\fa})=d$. Let $I$ be a
two-generated ideal of $A$ such that $IA/\fa$ is
$A/\fa$-regular. Then, for every $n\geq d$,
\begin{eqnarray*}
I^{n}\cap \fa =I^{n-d}(I^{d}\cap \fa) +\fm I^{n}\cap \fa .
\end{eqnarray*}
\end{Theorem}

\demo Let $I$ be a two-generated ideal of $A$ such that $IA/\fa$ is
$A/\fa$-regular.  By Proposition~\ref{slrt},
$s_{\mathfrak{m}}(\fa,A;I)\leq {\rm rt}_{\mathfrak{m}}(I;A/\fa)$. By
Remark~\ref{coef-mod}, ${\rm rt}_{\mathfrak{m}}(I;A/\fa)={\rm
rt}_{\mathfrak{m}/\mathfrak{a}}(IA/\fa)$, which is $d$ or less by
Proposition~\ref{rtm}.  \qed

\medskip

The next example, taken from Wang in \cite{wang}, shows that even this
weaker uniform Artin-Rees property of Theorem~\ref{ideals2} is not
true anymore for the set of three-generated ideals. It also shows that
if in Theorem~\ref{ideals} one changes the set of ideals having
principal reductions for the set of ideals having reductions generated
by regular sequences of length two, then there may not exist a uniform
Artin-Rees (modulo $\fm$) number.

\begin{Example}\label{wange} {\rm There exist $(A,\fm)$, a
three-dimensional noetherian local ring with infinite residue field,
$\fa$, a prime ideal of $A$ such that $A/\fa$ has finite integral
closure $\overline{A/\fa}$, and $\{I_{n}\}_{n}$, a family of
three-generated ideals of $A$ such that $I_{n}A/\fa$ is
$A/\fa$-regular, but, for every $n\geq 1$,
\begin{eqnarray*}
I_{n}^{n}\cap \fa\varsupsetneq I_{n}(I_{n}^{n-1}\cap \fa) +\fm
I_{n}^{n}\cap \fa . 
\end{eqnarray*}}
\end{Example}

\demo Take $(A,\fm)$, a three-dimensional regular local ring with
infinite residue field, $\fm=(x,y,z)$, the maximal ideal generated by
a regular system of parameters $x,y,z$, and $\fa=(z)$. Let
$I_{n}=(x^{n},y^{n},x^{n-1}y+z^{n})$. Since
$x^{n},y^{n},x^{n-1}y+z^{n}$ is a regular sequence of $A$, the
relation type of $I_{n}$ is ${\rm rt}(I_{n})=1$. It is not difficult
to prove that the relation type of $I_{n}A/\fa$ and that the relation
type of its fiber cone are given by ${\rm rt}(I;A/\fa)={\rm
rt}_{\mathfrak{m}}(I;A/\fa)=n$. Then, by Proposition~\ref{slrt},
$s_{\mathfrak{m}}(\fa,A;I_{n})={\rm
rt}_{\mathfrak{m}}(I_{n};A/\fa)=n$. Remark that $(x^{n},y^{n})A/\fa$
is a reduction of $I_{n}A/\fa$ generated by a regular sequence of
length two. \qed

\medskip

{\em Acknowledgement.} The second author is partially supported by the
MTM2004-01850 spanish grant.

\medskip 

\renewcommand{\baselinestretch}{1}

\medskip

\qquad
\parbox[t]{6cm}{\small

\noindent Jos\'e M. Giral 

\noindent Departament d'\`Algebra i Geometria

\noindent Universitat de Barcelona

\noindent Gran Via de les Corts Catalanes 585 

\noindent 08007 Barcelona, Spain

\noindent giral@ub.edu 
}
\qquad
\parbox[t]{7cm}{\small

\noindent Francesc Planas-Vilanova

\noindent Departament de Matem\`atica Aplicada 1 

\noindent Universitat Polit\`ecnica de Catalunya

\noindent Diagonal 647, ETSEIB

\noindent 08028 Barcelona, Spain

\noindent francesc.planas@upc.edu

}
{\small

\begin{thebibliography}{cc}
\bibitem{aht}{I. Aberbach, C. Huneke, N.V. Trung, Reduction numbers,
Brian\c{c}on-Skoda theorems and the depth of Rees rings, Compositio
Math. {\bf 97} (1995), no. 3, 403--434.}
\bibitem{akizuki}{Y. Akizuki, Einige Bemerkungen \"uber prim\"are
Integrit\"atsbereiche mit Teilerkettensatz, Proc. Phys.-Math. Soc.
Japan, {\bf 17} (1935), 327--336.}
\bibitem{am}{M.F. Atiyah, I.G. Macdonald, Introduction to
commutative algebra. Addison-Wesley Publishing Co., Reading,
Mass.-London-Don Mills, Ont. 1969.}
\bibitem{bourbaki}{ N. Bourbaki, El\'ements de math\'ematique.
Alg\`ebre. Chapitres 1 \`a 3. (French) Hermann, Paris 1970}
\bibitem{bs}{M.P. Brodmann, R.Y. Sharp, Local Cohomology. An algebraic
introduction with geometric applications, Cambridge Studies in
Advanced Mathematics, 60. Cambridge University Press, Cambridge,
1998.}
\bibitem{bh}{W. Bruns, J. Herzog, Cohen-Macaulay rings, Cambridge
Studies in Advanced Mathematics, 39. Cambridge University Press,
Cambridge, 1993.}
\bibitem{cpv}{A. Corso, C. Polini, W.V. Vasconcelos, Multiplicity of
the special fiber of blowups, Math. Proc. Cambridge Philos. Soc. {\bf
140} (2006), no. 2, 207--219.}
\bibitem{cz}{T. Cortadellas, S. Zarzuela, On the depth of the fiber
cone of filtrations, J. Algebra {\bf 198} (1997), no. 2, 428--445.}
\bibitem{costa}{D. Costa, Sequences of linear type, J. Algebra {\bf
94} (1985), 256-263.}
\bibitem{dgh1}{M. D'Anna, A. Guerreri, W. Heinzer, Invariants of ideals
having principal reductions, Comm. in Algebra {\bf 29}, no. 2, (2001),
889-906.}
\bibitem{dgh2}{M. D'Anna, A. Guerreri, W. Heinzer, Ideals having a
one-dimensional fiber cone, Ideal theoretic methods in commutative
algebra (Columbia, MO, 1999), 155--170, Lecture Notes in Pure and
Appl. Math., 220, Dekker, New York, 2001.}
\bibitem{dkv}{C. D'Cruz, V. Kodiyalam, J.K. Verma, Bounds on the
$a$-invariant and reduction numbers of ideals, J. Algebra {\bf 274}
(2004), no. 2, 594--601.}
\bibitem{do}{A.J. Duncan, L. O'Carroll, A full uniform
Artin-Rees theorem, J. reine angew. Math. {\bf 394} (1989),
203-207.}
\bibitem{eh}{D. Eisenbud, M. Hochster, A Nullstellensatz with
nilpotents and Zariski's main lemma on holomorphic functions,
J. Algebra {\bf 58} (1979), 157-161.}
\bibitem{hw}{W. Hassler, R.Wiegand, Direct sum cancellation for
modules over one-dimensional rings, J. Algebra {\bf 283} (2005),
93-124.}
\bibitem{hoa}{L.T. Hoa, Reduction numbers of equimultiple ideals.
J. Pure Appl. Algebra {\bf 109} (1996), no. 2, 111--126.}
\bibitem{huckaba1}{S. Huckaba, Reduction numbers for ideals of
analytic spread one, J. Algebra {\bf 108} (1987), no. 2, 503--512.}
\bibitem{huckaba2}{S. Huckaba, Reduction numbers for ideals of higher
analytic spread, Math. Proc. Cambridge Philos. Soc. {\bf 102} (1987),
49-57.}
\bibitem{huckaba3}{S. Huckaba, On complete $d$-sequences and the
defining ideals of Rees algebras, Math. Proc. Cambridge Philos. Soc.
{\bf 106} (1989), no. 3, 445--458.}
\bibitem{huneke}{C. Huneke, Uniform bounds in Noetherian
rings, Invent. Math. {\bf 107} (1992), no. 1, 203--223.}
\bibitem{hs}{C. Huneke, I. Swanson, Integral closure of ideals, rings,
and modules, London Mathematical Society Lecture Note Series,
336. Cambridge University Press, Cambridge, 2006.}
\bibitem{matsumura}{H. Matsumura, Commutative Ring Theory, Cambridge
Studies in Advanced Mathematics, Cambridge University Press, 1986.}
\bibitem{nagata}{M. Nagata, Local rings. Interscience Tracts in Pure
and Applied Mathematics, No. 13 Interscience Publishers a division of
John Wiley \& Sons, New York-London 1962.}
\bibitem{nr}{D.G. Northcott, D. Rees, Reductions of ideals in local
rings, Proc. Cambridge Philos. Soc. {\bf 50} (1954), 145--158.}
\bibitem{ocarroll1}{ L. O'Carroll, A uniform Artin-Rees theorem and
Zariski's main lemma on holomorphic functions, Invent. Math. {\bf 90}
(1987), 647--652.}
\bibitem{ocarroll2}{ L. O'Carroll, A note on Artin-Rees
numbers, Bull. London Math. Soc. {\bf 23} (1991), 209-212.}
\bibitem{planas}{F. Planas-Vilanova, The strong uniform Artin-Rees
property in codimension one, J. reine angew. Math. {\bf 527} (2000),
185-201.}
\bibitem{reid}{M. Reid, Undergraduate commutative algebra. Cambridge
University Press, Cambridge, 1995.}
\bibitem{rossi}{M.E. Rossi, A bound on the reduction number of a
primary ideal, Proc. Amer. Math. Soc. {\bf 128} (2000), no. 5,
1325--1332.}
\bibitem{sv}{J. Sally, W.V. Vasconcelos, Stable rings.  J. Pure
Appl. Algebra {\bf 4} (1974), 319--336.}
\bibitem{schenzel}{P. Schenzel, Castelnuovo's index of regularity and
reduction numbers.  Topics in algebra, Part 2 (Warsaw, 1988),
201--208, Banach Center Publ., 26, Part 2, PWN, Warsaw, 1990.}
\bibitem{trung1}{N.V. Trung, Reduction exponent and degree bound for
the defining equations of graded rings, Proc. Amer. Math. Soc. {\bf
101}, no. 2, (1987), 229-236.}
\bibitem{trung2}{N.V. Trung, The Castelnuovo regularity of
the Rees algebra and the associated graded ring,
Trans. Amer. Math. Soc. {\bf 350} (1998), 2813-2832.}
\bibitem{vasconcelos1}{W.V. Vasconcelos, Reduction numbers of ideals,
J. Algebra {\bf 216} (1999), no. 2, 652--664.}
\bibitem{vasconcelos2}{W.V. Vasconcelos, Multiplicities and reduction
numbers, Compositio Math. {\bf 139} (2003), no. 3, 361--379.}
\bibitem{vasconcelos3}{W.V. Vasconcelos, Integral closure. Rees
algebras, multiplicities, algorithms. Springer Monographs in
Mathematics. Springer-Verlag, Berlin, 2005.}
\bibitem{wang}{H.J. Wang, Some uniform properties of 2-dimensional
local rings, J. Algebra {\bf 188} (1997), 1-15.}
\end{thebibliography}}
\end{document}